
\mag=1000 \hsize=6.8 true in \vsize=8.7 true in
 \baselineskip=15pt
\vglue 1.5 cm

       


 \def\qed{$\rlap{$\sqcap$}\sqcup$}


\font\tengothic=eufm10 \font\sevengothic=eufm7
\newfam\gothicfam
     \textfont\gothicfam=\tengothic
     \scriptfont\gothicfam=\sevengothic
\def\goth#1{{\fam\gothicfam #1}}


  \font\tenmsb=msbm10              \font\sevenmsb=msbm7
\newfam\msbfam
     \textfont\msbfam=\tenmsb
     \scriptfont\msbfam=\sevenmsb
\def\Bbb#1{{\fam\msbfam #1}}



\def\move-in{\parshape=1.75true in 5true in}


\hyphenation {Castel-nuovo}


\def\PP#1{{\Bbb P}^{#1}}

\def\PPP {{\Bbb P}}
\def\ref#1{[{\bf #1}]}
\def\prim{^\prime }

\def\Proof{\noindent {\it Proof:}}
\def\Prop#1{\noindent {\bf Proposition #1:}}


\def\Coroll#1{\noindent {\bf Corollary #1:}}



\vskip 1cm \centerline{\bf Higher Secant Varieties of 
Segre-Veronese varieties.}
\bigskip

\centerline{\it M.V.Catalisano, A.V.Geramita, A.Gimigliano.}
\bigskip
\bigskip
\noindent {\bf 0. Introduction.}
\par
\medskip
\bigskip
The problem of determining the dimensions of the higher secant
varieties of the classically studied projective varieties (and to
describe the defective ones) is a problem with a long and
interesting history.

In the case of the Segre varieties there is much interest in this
question, and not only among geometers.  In fact, this particular
problem is strongly connected to questions in representation
theory, coding theory, algebraic complexity theory (see our paper
\ref{CCG2} for some recent results as well as a summary of known
results, and also \ref{BCS}) and, surprisingly enough, also in
algebraic statistics (e.g. see \ref {GHKM} and \ref {GSS}).

We address a slight generalization of this problem here; more precisely we will study the
higher secant varieties of $${\Bbb X} = \PP {n_1}\times ...\times
\PP {n_t}= \PP {\bf n}, \ {\bf n} =(n_1,...,n_t)$$  embedded in
the projective space $\PP {N}$ ($N= \Pi {a_i+n_i\choose n_i}-1$)
by the morphism ${\nu }_{\bf n, a}$ given by $ {\cal O} _{\PP {n}}({\bf a})$,
where ${\bf a}=(a_1,...,a_t)$ ($a_i$ positive integers).  We
denote the embedded variety ${\nu }_{\bf n, a}({\Bbb X})$ by $V_{{\bf n},{\bf a}}$,
and call it a {\it Segre-Veronese variety} and the embedding a
{\it Segre-Veronese embedding} (see \ref {BuM}).

In Section 1 we recall some classical results by Terracini
regarding  secant varieties and we also introduce one of the
fundamental observations ({\it Theorem 1.1}) which allows us to
convert certain questions about ideals of varieties in
multiprojective space to questions about ideals in standard
polynomial rings. 

In Section 2 we concentrate on $t = 2,3$ and let ${\bf a}$ be
arbitrary.  In {\it Theorem 2.1} we give the dimensions of all the
higher secant varieties for $V_{{\bf n}, {\bf a}}$ (where ${\bf n}
= (1,1)$, and $\bf a$ is arbitrary).  For ${\bf n} = (k,n), \ {\bf
a} = (1,k+1)$ we find that $V_{{\bf n},{\bf a}}$ has no deficient
higher secant varieties ({\it Proposition 2.3}) and this gives an
interesting conclusion about the Grassmann defectivity for the
$(k+1)$-Veronese embedding of $\PP {n}$ ($k\geq 2, \ n\geq 2$).
We also state our theorem (the proof will appear elsewhere) which
gives the dimensions of all the higher secant varieties to
$V_{(1,1,1), (a,b,c)}$ for any positive integers $a,\ b$ and $c$
({\it Theorem 2.5}).

Section 3 is dedicated to results on regularity of secant varieties of 
Segre Veronese varieties which can be deduced from results for the Segre varieties and by studying the multigraded Hilbert function of a scheme of 2-fat points in $\PP {\bf n}$. Also we give several examples of defective and Grassmann defective Segre-Veronese varieties.

Finally, in  Section 4 we describe a way of thinking about the points of
Segre-Veronese varieties as (partially symmetric) tensors.

Our method is essentially this (see \S 1): we use Terracini's
Lemma (as in \ref {CGG2}, \ref {CGG4}) to translate the problem of
determining the dimensions of higher secant varieties into that of
calculating the value, at $(a_1,...,a_t)$, of the Hilbert function
of generic sets of 2-fat points in $\PP {\bf n}$. Then we show, by
passing to an affine chart in $\PP {\bf n}$ and then homogenizing
in order to pass to $\PP n$, $n=n_1+...+n_t$, that this last
calculation amounts to computing the Hilbert function of a very
particular subscheme of $\PP n$.  Finally, we study the
postulation of these special subschemes of $\PP {n}$.


We wish to warmly thank Monica Id\`a, Luca Chiantini and Ciro
Ciliberto for many interesting conversations about the questions
considered in this paper.

\bigskip\noindent {\bf 1. Preliminaries, the Multiprojective-Affine-Projective Method.}

\bigskip Let us recall the notion of higher secant varieties.

\medskip\noindent {\bf Definition 1:} Let $X\subseteq \PP N$ be a closed
irreducible projective variety of dimension $n$. The $s^{th}$ {\it
higher secant variety} of $X$, denoted $X^s$, is the closure of
the union of all linear spaces spanned by $s$ independent points
of $X$.

Recall that, for $X$ as above, there is an inequality involving
the dimension of $X^s$.  Namely, $$ \dim X^s \leq {\rm min} \{N,
sn+s-1\}  \eqno (1) $$ and one ``expects" the inequality should,
in general, be an equality.

When $X^s$ does not have the expected dimension, $X$ is said to be
$(s-1)$-{\it defective}, and the positive integer $$ \delta
_{s-1}(X) := {\rm min} \{N, sn+s-1\}-\dim X^s $$ is called the
$(s-1)${\it-defect} of $X$.  Probably the most well known
defective variety is the Veronese surface, $X$ in $\PP 5$ for
which $\delta_1(X) = 1$.

As a generalization of the higher secant varieties of a variety,
one can also consider the following varieties.

\medskip\noindent{\bf Definition 2:}\ Let $X\subseteq \PP N$ be a closed
irreducible projective variety of dimension $n$. The {\it
(k,s-1)-Grassmann secant variety,} denoted $Sec_{k,s-1}(X)$, is
the Zariski closure (in the Grassmaniann of $k$-dimensional linear
subspaces of $\PP N$) of the set $$ \{l\in {\Bbb G}(k,N)\ \vert \
l\ lies\ in\ the\ span\ of\ s\ independent\ points\ of\ X\}. $$ In
case $k=0$ we get $X^s = Sec_{0,s-1}(X)$.

As a generalization of the analogous result for the higher secant
varieties, one always has $$ \dim Sec_{k,s-1}(X) \leq {\rm min}\{
sn+(k+1)(s-k-1), (k+1)(N-k)\}, $$ with equality being what is
generally ``expected".

When  $Sec_{k,s-1}(X)$ does not have the expected dimension then
we say that $X$ is $ (k,s-1)${\it-defective} and in this case we
define the $(k,s-1)$-defect of $X$ as the number: $$ \delta
_{k,s-1}(X)= {\rm min}\{ sn+(k+1)(s-k-1), (k+1)(N-k)\}-\dim
Sec_{k,s-1}(X). $$ (For general information on these defectivities
see {\ref {ChCo}} and  \ref {DF}.)

\medskip In this note we study the defectivities of $V_{{\bf n},{\bf a}} $.  For convenience, and when no doubts can arise about the variety we are considering, we will just write $V$ for $V_{{\bf n},{\bf a}} $.

In his paper \ref {Te2}, Terracini gives a link between these two
kinds of defectivity for a variety $X$ as above (see \ref {DF} for
a modern proof):

\medskip \Prop {1.0} {\it (Terracini)} {\it Let $X\subset \PP N$ be an irreducible non-degenerate
projective variety of dimension n. Let $\sigma : X \times \PP k
\rightarrow \PP {(k+1)(N+1)-1}$ be the (usual) Segre embedding. Then X
is $(k,s-1)$-defective with defect $\delta _{k,s-1}(X)=\delta$ if
and only if $\sigma (X \times \PP k)$ is $(s-1)$-defective with
$(s-1)$-defect $\delta_{s-1}(X \times \PP k)=\delta$.}

\bigskip A classical result about higher secant varieties is Terracini's Lemma (see \ref {Te}, \ref{CGG2}):

\medskip\noindent {\bf Terracini's Lemma:} {\it Let} $(X,{\cal L})$ {\it be a
polarized, integral, non-singular scheme; if} ${\cal L}$ {\it
embeds X into} $\PP N$ {\it , then:} $$ T_P(X^s) =
<T_{P_1}(X),...,T_{P_s}(X)>, $$ {\it where} $P_1,...,P_s$ {\it are
s generic points on X, and P is a generic point of $<P_1,...,P_s>$
{\it (the linear span of $P_1, \ldots , P_s$); here }
$T_{P_i}(X)$ {\it is the projectivized tangent space of X in} $\PP
N$.} {\hfill \qed}

\medskip  Let $Z\subset X$ be a scheme of $s$  generic 2-{\it fat points}, that is a scheme defined by the ideal sheaf ${\cal I}_{Z} = {\cal I}^2_{P_1}\cap ...\cap {\cal I}^2_{P_s}\subset {\cal O}_X$, where $P_1,...,P_s$  are $s$ generic points.  Since there is a bijection between hyperplanes of the space $\PP N$ containing the subspace $<T_{P_1}(X),...,T_{P_s}(X)>$ and the elements of $H^0(X,{\cal I}_{Z}({\cal L}))$, we have:

\medskip \Coroll {} {\it Let X,} ${\cal L}$, $Z$ {\it , be as above; then}
$$ \dim X^s = \dim <T_{P_1}(X),...,T_{P_s}(X)> = N - \dim
H^0(X,{\cal I}_{Z}({\cal L})) . $$ {\hfill \qed}

\medskip\noindent  Let $X=\PP {n_1}\times ... \times \PP {n_t}$ and let $V_{{\bf n},{\bf a}}= V\subset \PP N$ be the
embedding of $X$ given by ${\cal L}={\cal O}_X(a_1,...,a_t)$.  By applying the corollary above to our case (i.e.
$V=V_{{\bf n},{\bf a}}$), we get
$$ \dim V^s = H(Z,{\bf a})-1, \eqno (2) $$ where $Z\subset  {\Bbb
P}^{n_1}\times ...\times \PP {n_t}$ is a set of $s$ generic 2-fat
points, and where $\forall {\bf j}\in {\Bbb N}^t$, $H(Z,{\bf j})$
is the Hilbert function of $Z$ i.e. $$ H(Z,{\bf j}) = \dim R_{\bf
j} - \dim H^0({\Bbb P}^{n_1}\times ...\times \PP {n_t},{\cal
I}_{Z}({\bf j})), $$ where $R=k[x_{0,1},...,x_{n_1,1},\ ...\
,x_{0,t},...,x_{n_t,t}]$ is the multi-graded homogeneous
coordinate ring of $\PP {n_1}\times ...\times \PP {n_t}$.

Now let $n=n_1+...+n_t$ and consider the birational map $$ g: \PP
{n_1}\times ...\times \PP {n_t} ---\rightarrow {\Bbb A}^n, $$
where: $$ ((x_{0,1},...,x_{n_1,1}),...,(x_{0,t},...,x_{n_t,t}))
\longmapsto ({x_{1,1}\over x_{0,1}},{x_{2,1}\over
x_{0,1}},...,{x_{n_1,1}\over x_{0,1}}; {x_{1,2}\over
x_{0,2}},...,{x_{n_2,2}\over x_{0,2}};...; {x_{1,t}\over
x_{0,t}},...,{x_{n_t,t}\over x_{0,t}}) \ . $$ This map is defined
in the open subset of $\PP {n_1}\times ...\times \PP {n_t}$  given
by $\{x_{0,1}x_{0,2}...x_{0,t}\neq 0\}$.

Let $S=k[z_0,\ z_{1,1},...,z_{n_1,1},\ z_{1,2},...,z_{n_2,2},\
...\ , z_{1,t}, ... \ ,z_{n_t,t}]$ be the coordinate ring of $\PP
n$ and consider the embedding ${\Bbb A}^n \rightarrow \PP n$ whose
image is the chart ${\Bbb A}^n_0=\{z_0=1\}$. By composing the two
maps above we get: $$ f: \PP {n_1}\times ...\times \PP {n_t}
---\rightarrow {\Bbb P}^n, $$ with $$
((x_{0,1},...,x_{n_1,1}),...,(x_{0,t},...,x_{n_t,t})) \longmapsto
(1,{x_{1,1}\over x_{0,1}},...,{x_{n_1,1}\over
x_{0,1}};{x_{1,2}\over x_{0,2}},...,{x_{n_2,2}\over x_{0,2}};...;
{x_{1,t}\over x_{0,t}},...,{x_{n_t,t}\over x_{0,t}}) $$ $$ =
(x_{0,1}x_{0,2}...x_{0,t},\ x_{1,1}x_{0,2}...x_{0,t},\
x_{0,1}x_{1,2}...x_{0,t},\ ...,\ x_{0,1}...x_{0,t-1}x_{n_t,t}). $$

Let $Z\subset \PP {n_1}\times ...\times \PP {n_t}$ be a
zero-dimensional scheme which is contained in the affine chart
$\{x_{0,1}x_{0,2}...x_{0,t}\neq 0\}$ and let $Z\prim = f(Z)$. We
want to construct a scheme $W\subset  \PP n$ such that $\dim
(I_W)_a=\dim (I_Z)_{(a_1,...,a_t)}$, where $a=a_1+...+a_t$.

Let us recall that the coordinate ring of $\PP n$ is
$S=k[z_0,\ z_{1,1},...,z_{n_1,1},\ z_{1,2},...,z_{n_2,2},\
...\ , z_{1,t}, ... \ ,z_{n_t,t}]$, and let
$Q_0,Q_{1,1},...,Q_{n_1,1},Q_{1,2},...,Q_{n_2,2},...,Q_{n_t,t}$ be
the coordinate points of $\PP {n}$. Consider the linear
subspace $\Pi_i\cong \PP {n_i-1}  \subset \PP {n}$, where
$\Pi_i=<Q_{1,i},...,Q_{n_i,i}>$.  The defining ideal of $\Pi_i$
is: $$ I_{\Pi_i} = (z_0, \ z_{1,1}, \ldots , z_{n_1,1};\ \ldots ;
\hat z_{1,i},\ldots ,\hat z_{n_i,i}; \ \cdots ;\ z_{1,t}, \cdots ,
z_{n_t,t}) \ . $$ Let $W_i$ be the subscheme of $\PP {n}$ denoted
by $(a-a_i)\Pi_i$, i.e. the scheme defined by the ideal
$I_{\Pi_i}^{a-a_i}$. Notice that $W_i\cap W_j = \emptyset$ for
$i\neq j$.

\medskip\noindent {\bf Theorem 1.1:} {\it Let $Z,\ Z\prim$ be as above and let
$W = Z\prim + W_1 + ... + W_t\subset \PP n$.  Then we have: $$
\dim (I_W)_a=\dim (I_Z)_{(a_1,...,a_t)} $$ where $a=a_1+...+a_t$.}

\medskip\noindent  {\it Proof:} First note that
$$ R_{(a_1, \ldots, a_t)} = \langle
(x_{0,1}^{a_1-s_1}N_1)(x_{0,2}^{a_2-s_2}N_2)\cdots
(x_{0,t}^{a_t-s_t}N_t) \rangle $$ where every $N_i$ varies among all  monomials in
$(x_{1,i}, \ldots , x_{n_i,i})_{s_i}$, for all $s_i \leq a_i$.

By dehomogenizing (via $f$ above) and then substituting $z_{i,j}$
for $(x_{i,j}/x_{0,j})$, and finally homogenizing with respect to
$z_0$, we see that $$ R_{(a_1, \ldots , a_t)} \simeq \langle
z_0^{a-s_1-\cdots - s_t}M_1M_2\ldots M_t \rangle $$ where every
$M_i$ varies among all monomials in $(z_{1,i}, \ldots , z_{n_i,i})_{s_i}$.

\medskip\move-in\noindent{\it Claim:}\ 
$$ (I_{W_1 +...+ W_t })_a= (I_{W_1} \cap ...\cap I_{W_t} )_a=
\langle  z_0^{a-s_1-\cdots - s_t}M_1M_2\ldots M_t \rangle $$
where every $M_i$ varies  among all monomials in $(z_{1,i}, \ldots , z_{n_i,i})_{s_i}$, 
for all $s_i \leq a_i$.

\medskip
\move-in\noindent{\it Proof of Claim:} \ $\subseteq$: \ Since both
vector spaces are generated by monomials, it is enough to show
that the monomials of the left hand side of the equality are
contained in the right hand side of the equality.

\move-in Consider $M = z_0^{a-s_1-\cdots - s_t}M_1M_2\ldots M_t$
(as above).  We now show that this monomial is in $I_{W_i}$ (for
each $i$).  Notice that $M_j \in I_{\Pi_i}^{s_i}$ (for $j \neq i$)
and that $z_0^{a-s_1-\cdots - s_t} \in I_{\Pi_i}^{a-s_1-\cdots -
s_t}$.  Thus, $M \in I_{\Pi_i}^{(a-s_1-\cdots - s_t) + (s_1 +
\cdots + \hat{s_i} + \cdots + s_t)} = I_{\Pi_i}^{a-s_i}$.  Since
$s_i \leq a_i$ we have $a-a_i \leq a-s_i$ and so $M \in
I_{\Pi_i}^{a-a_i}$ as well, and that is what we wanted to show.

\move-in\medskip\noindent $\supseteq$\ : To prove this inclusion,
consider an arbitrary monomial $M \in S_a$.  Such an $M$ can be
written $M = z_0^{\alpha _0}M_1\cdots M_t$ where $M_i \in
(z_{1,i}, \ldots , z_{n_i,i})$ is a monomial of degree $\alpha
_i$.

\move-in  Now, $M \in (I_{W_1 +...+ W_t})_a$ means $M \in (I_{W_i})_a$ for each
$i$, hence   $$ \alpha _0 + \alpha _1 + \cdots + \hat{\alpha _i} +
\cdots + \alpha _t \geq a - a_i $$ for $i = 1, \ldots , t$.  Since
$ \alpha _0 + \alpha _1 + \cdots + \alpha _t = a$,
then $a-\alpha _i \geq a - a_i$ for each $i$, and so $\alpha
_i \leq a_i$ for each $i$.  That finishes the proof of the Claim.

\medskip Now, since $Z$ and $Z\prim$ are isomorphic 
($f$ is an isomorphism between the two affine charts $\{z_0\ne
q 0\}$ 
and $\{x_{0,1}x_{0,2}...x_{0,t}\neq 0\}$), it
immediately follows (via the two different dehomogeneizations) that $(I_Z)_{(a_1,...,a_t)}\cong (I_W)_a$.  
\par \hfill  \qed

\medskip When $Z$ is given by $s$ generic 2-fat points, we have the obvious corollary:

\medskip\Coroll {1.2} {\it Let $Z\subset \PP {n_1}\times ...\times \PP {n_t}$ be a
generic set of s 2-fat points, let $W\subset \PP n$ be as in
Theorem 1.1, then we have:} $$ \dim  V^s = H(Z,(a_1,...,a_t))-1 =
N - \dim (I_W)_a. $$

\bigskip 

\bigskip
\noindent {\bf 2. On Segre-Veronese with two or three factors.}
\par
\medskip
First we consider the case $\PP 1\times \PP 1$, i.e. $t=2$, $n_1=n_2=1$,
for all ${\bf a}=(a_1,a_2)\in {\Bbb N}^2$.
We get that the $\Pi_i$'s are points; let $\Pi_1=
A_1=(0,1,0)$, $\Pi_2=A-2=(0,0,1)$, in $\PP 2$, and let $Z =2R_1 + ... +
2R_s \subset  \PP {1}\times \PP {1}$ be a set of $s$ generic 2-fat
points. We may assume that $R_i=((1,\alpha_i),(1,\beta_i))$, so
that $f: \PP 1\times \PP 1 - - - \rightarrow \PP 2$ is such that:
$$f(R_i)=P_i =(1,\alpha_i,\beta_i)\in \PP 2,$$ and
$A_1,A_2,P_1,...,P_s$  will be generic points of $\PP 2$. Let
$$W=a_2A_1+a_1A_2+2P_1 + ... + 2P_s \subset \PP 2$$ be the scheme defined
by the ideal sheaf ${\cal I}_{W} = {\cal I}^{a_2}_{A_1} \cap 
{\cal I}^{a_1}_{A_2} \cap {\cal I}^2_{P_1}\cap ...\cap {\cal I}^2_{P_s}$. By
Theorem 1.1, the data of the Hilbert function of $Z$ in bidegree
$(a_1,a_2)$ is equivalent to the data of the Hilbert function of $W$
in degree $(a_1+a_2)$, in fact, from $(I_W)_{a_1+a_2}\cong (I_Z)_{(a_1,a_2)}$,
we easily get: $$ H(W,{(a_1+a_2)})=H(Z,{(a_1,a_2)}) + \deg(a_2A_1+a_1A_2)=
H(Z,{(a_1,a_2)}) + {a_1+1 \choose 2} + {a_2+1 \choose 2} $$ since $$ \dim
(I_Z)_{(a_1,a_2)} = (a_1+1)(a_2+1)-H(Z,{(a_1,a_2)}) $$ $$ \dim (I_W)_{(a_1+a_2)})
= {a_1+a_2+2 \choose 2} - H(W,{(a_1,a_2)}). $$ Now let
$V=V_{{\bf 1},{\bf a}}=\nu _{{\bf 1},{\bf a}}(\PP 1\times \PP 1)\subset \PP
{a_1a_2+a_1+a_2}$ be the Segre-Veronese embedding; by Corollary 1.3 we
get: $$ \dim V^s = H(Z,{(a_1,a_2)})-1 =a_1a_2+a_1+a_2 -\dim (I_W)_{(a_1+a_2)}= $$
$$ = H(W,{(a_1+a_2)})-1-{a_1+1 \choose 2}-{a_2+1 \choose 2} $$
\medskip
Hence, in order to compute $\dim V^s$, we should study a scheme of
generic fat points $W$ in $\PP {2}$; without loss of genericity,
we may suppose $a_1\geq a_2$.
\medskip
\medskip
\noindent {\bf Theorem 2.1:} {\it Let $V=V_{{\bf 1},{\bf a}}=\nu
_{{\bf 1},{\bf a}}(\PP 1\times \PP 1)$, then $V^s$ has the expected
dimension, except for } $$ a_1=2d, \   a_2=2,\ d\geq 1,\  {\rm and} \
\ s=a_2+1 \eqno (\dag) $$ {\it In this case $V^s$ is defective, and
$\dim V^s =3s-2$ (its defectivness is $1$)}.
\par
\medskip
This theorem could be proved by methods similar to those used in 
\ref {CGG3}, \ref {CGG4}, but we omit the lenghty and tedious proof here. Notice that in \ref {ChCi} a total classification of all the surfaces
with some defective secant variety can be found and in \ref {La} the case of rational scroll is treated, which covers also the $\PP 1\times \PP 1$ case. We wish to thank Monica Id\`a for showing us a direct proof of Theorem 2.1 (\ref {Id}) which uses the "differential Horace" method.
\par
\bigskip
Theorem 2.1 yields (by interpreting things via Theorem 1.1) that
for fat points in $\PP 2$, we have:
\par
\medskip
\noindent {\bf Remark 2.2:} Let $a_1,a_2,s$ be positive integers, with
$a_1\geq a_2$. Let $W=a_2A_1+a_1A_2+2P_1+...+2P_s\subset \PP 2$. Then $$
H(W,a_1+a_2) = min \{{a_1+a_2+2\choose 2}, {a_1+1\choose 2}+{a_2+1\choose
2}+3s\}, \eqno (*) $$ except when $a_1=2d$, $a_2=2$ and $s=2d+1$; in
this case $H(W,a_1+a_2)$ is $1$ less than expected.
\par
\medskip
Notice that in the exceptional case it is easy to check what is
the geometrical situation: there is an unique (rational) curve $C$
through $dA_1+A_2+P_1+...+P_{2d+1}$, and $2C$ gives an ``unexpected
element" of $(I_W)_{2d+2}$.
\par
\bigskip
\bigskip
Now let us consider another case; namely the products $\PP r\times
\PP k$; we have the following:
\par
\bigskip
\Prop {2.3} {\it Let $r,k\geq 1$, and $V\subset \PP N$ be the
$(k+1,1)$ Segre-Veronese embedding of $\PP r\times \PP k$. Then
for any $s\geq 1$, $V^s$ has the expected dimension.}
\par
\medskip
\Proof \  Notice that, since $N={r+k+1\choose
r}(k+1)-1$,  for $s={r+k\choose k}$ we get
that the expected dimension of $V^s$ is exactly $N$, hence the
statement will follow if we prove that $V^{r+k\choose k}=\PP N$.
\par
Consider the scheme $W=\Pi _1+ (k+1)\Pi_2+ 2P_1+ ...+
2P_s\subset \PP {r+k}$, where $\Pi_1\cong \PP {r-1}$ and
$\Pi_2\cong \PP {k-1}$ are linear spaces and $s={r+k\choose k}$;
then, by Corollary 1.2, we get that $$\dim V^s = N - h^0(\PP
{r+k},{\cal I}_W(k+2)). $$ From what we have seen before, we will
be done if $\dim (I_W)_{k+2} = {0}$.
\par
We will proceed by double induction on $k$ and $r$.
\par
Let us consider the case $k=1$ (any $r$) first. When $k=1$, $W$ is
the schematic union: $W=\Pi _1+ 2\Pi_2+ 2P_1+ ...+
2P_s\subset \PP {r+1}$; where $\Pi_2 $ is a point. It is
enough, in order to prove the case $k=1$, to show that $(I_W)_{3}
= \{0\}$ (here $s=r+1$).
\par
Let us work by induction on $r$; for $r=1$ we trivially have
$h^0(\PP {2},{\cal I}_W(3))=0$ (since $W$ is made of three 2-fat
points and one simple point). When $r>1$, consider the exact
sequence: $$ 0 \rightarrow {\cal I}_{W'}(2) \rightarrow {\cal
I}_W(3)\rightarrow {\cal I}_{W\cap H,H}(3)\rightarrow 0 $$ where
$H\subset \PP{r+1}$ is the hyperplane $H=<P_1,...,P_{r+1}>$ and
$W'=\Pi_1+2\Pi_2+P_1+...+P_{r+1}$.  We get $h^0(H,{\cal I}_{W\cap
H,H}(3))=0$ by induction, and $h^0(\PP {r+1},{\cal I}_{W'}(2))=0$
since any of its element should give a quadric cone with vertex in
$\Pi_2$, but, since $\Pi_1\cong \PP {r-1}$, the cone should split
into the hyperplane $<\Pi_1,\Pi_2>$, and another hyperplane containing
$\Pi_2,P_1,...,P_{r+1}$, which is impossible by their genericity.
\par
Hence the case $k=1$ is done.
\par
\medskip
Now let us consider the case $r=1$; here we have $\PP 1\times \PP
k\rightarrow \PP {N}$, and $W=\Pi_1+(k+1)\Pi
_2+2P_1+...+ 2P_{k+1}\subset \PP {k+1}$, where $\Pi
_1$ is a point  and we must show that $(I_W)_{k+2}=\{0\}$.
\par
The hyperplanes $H_i=<P_i,\Pi_2>$, $i=1,...,k$, are fixed components for the hypersurfaces given by the forms in 
$({I}_W)_{k+2}$, hence (by removing such fixed components), we get 
$$ \dim ({I}_W)_{k+2} = \dim ({I}_{\Pi_1+P_1+...+P_{k+1}})_{1} =0.
$$
So case $r=1$ is done.
\par
\medskip
Now we can consider $r,k\geq2$ and work by double induction on
them.
\par
\medskip
Let $H\subset \PP {r+k}$ be a hyperplane such that $\Pi_2\subset
H$ and $\Pi _1$ is not contained in $H$. Let $\Pi_1'= H\cap \Pi_1
\cong \PP {r-2}$. Specialize $P_1,...,P_{s'}$, $s'={r+k-1\choose
k}$ on $H$ and consider the exact sequence: $$ 0\rightarrow {\cal
I}_{Z}(k+1)\rightarrow {\cal I}_W(k+2) \rightarrow {\cal I}_{W\cap
H,H}(k+2)\rightarrow 0, $$ where $Z=\Pi_1+ k\Pi_2+ P_1+
...+ P_{s'}+ 2P_{s'+1}+ ...+ 2P_s\subset \PP {r+k}$
and $W\cap H = \Pi_1'+ (k+1)\Pi_2+ 2P_1+ ...+
2P_{s'}\vert _H$.
\par
We have $h^0({\cal I}_{W\cap H,H}(k+2))=0$ by induction on
$r$; so we will be done if $h^0({\cal I}_{Z}(k+1))=0$, since the
above sequence would yield $h^0({\cal I}_W(k+2))=0$.
\par
Let us consider now a hyperplane $H'\subset \PP {r+k}$ with
$\Pi_1\subset H'$ and $\Pi _2$ not contained in $H'$; let
$\Pi_2'=H'\cap \Pi_2\cong \PP {k-2}$, then specialize
$P_{s'+1},...,P_s$ on $H'$ and consider the exact sequence: $$
0\rightarrow {\cal I}_{Z'}(k)\rightarrow {\cal I}_Z(k+1)
\rightarrow {\cal I}_{Z\cap H',H'}(k+1)\rightarrow 0, $$ where
$Z'=k\Pi_2+ P_1+ ...+ P_s\subset \PP {r+k}$ and $Z\cap H'
= \Pi_1+ k\Pi_2'+ 2P_{s'+1}+ ...+ 2P_{s}\vert _{H'}$.
Notice that $s-s'={r+k\choose k}-{r-1+k\choose k}= {r+k-1\choose
k-1}$, so that $h^0({\cal I}_{Z\cap H',H'}(k+1))=0$ by induction
on $r$ and $k$.
\par
So we are only left to prove $h^0({\cal I}_{Z'}(k))=0$.  The
sections of ${\cal I}_{Z'}(k)$ correspond to degree $k$
hypersurfaces in $\PP {r+k}$ which, in order to contain $k\Pi_2$
have to be cones with $\Pi_2$ as vertex. Let $H''\cong \PP r$ be a
generic $r$-dimensional linear subspace of $\PP {r+k}$; then
$h^0({\cal I}_{Z'}(k))=h^0({\cal I}_{Z'',H''}(k))$, where
$Z''\subset H''$ is the projection of $Z'$ into $H''$ from
$\Pi_2$. We have $Z''=Q_1+ ...+ Q_{s'}+ Q_{s'+1}+
...+ Q_s$, where $Q_{s'+1}, ..., Q_s$ are generic in $H''$,
while $Q_1,...,Q_{s'}$ are contained in the linear space $H\cap
H''\cong \PP {r-1}$ (where they are generic). A hypersurface $F$
of degree $k$ in $H''$ cannot contain $Q_1,...,Q_{s'}$ without
containing all $H\cap H''$ ($F$ intersects $H\cap H''$ in
something of degree $k$, but the $Q_i$, $i=1,...,s'={r+k-1\choose
k}$ are generic in $H\cap H''$, so they are not contained in a
hypersurface of degree $k$ of $\PP {r-1}=H\cap H''$). So, if $F$
vanishes on $Q_1,...,Q_s$, we have that $F$ is the union of $H\cap
H''$ and a hypersurface $F'$ of degree $k-1$ vanishing on
$Q_{s'+1}, ..., Q_s$, but again this cannot happen because they
are ${r+k-1\choose k-1}$ generic points in $\PP r$, so no form of
degree $k-1$ vanishes at them. Thus $h^0({\cal I}_{Z'',H''}(k))=0$
and we are done.
\par \hfill  \qed
\par
\medskip
From this result we get immediately the following:
\par
\Coroll {2.4} {\it Let $r,k\geq 2$, and $V$ be the (k+1)-ple
(Veronese) embedding of $\PP r$. Then $V$ is not (Grassman)
$(k,s-1)$-defective, for any s.}
\par
\medskip
\Proof\ By  Proposition 1.0, this statement is equivalent to
Proposition 2.3.\par \hfill  \qed
\par
\bigskip
\bigskip
In the case $t=3$, i.e. $\PP 1\times \PP 1\times \PP1= (\PP 1)^3$,
the situation for all Segre-Veronese embeddings can be analyzed
with the same methods used above; if we consider the embedding
$V_{\bf a}$ of $(\PP 1)^3$ given by the forms of tridegree ${\bf
a}=(a_1,a_2,a_3)$, in order to compute the dimension of
$V_{\bf a}^s$ we will have to study (by Theorem 1.1) the scheme of
fat points $$ W_s=(a_2+a_3)A_1+(a_1+a_3)A_2+(a_1+a_2)A_3+2P_1+...+2P_s \subseteq \PP
3, $$ where $A_1$, $A_2$, $A_3$ are coordinate points.
\par
A complete description of what happens is given by the following
result:
\par
\bigskip
\noindent {\bf Theorem 2.5:} {\it Let $a_1\geq a_2 \geq a_3\geq 1$,
$\alpha \in {\Bbb N}$ and $V=V_{\bf a}$ be a Segre-Veronese
embedding of $\PP 1\times \PP 1\times \PP 1$. Then $V^s$ has the
expected dimension, except for:
\par
$$ (a_1,a_2,a_3) = (2,2,2),\qquad {\rm and}\qquad s=7;$$ $$ (a_1,a_2,a_3) =
(2\alpha,1,1),\qquad {\rm and}\qquad s=2\alpha + 1.$$ In these
cases $V^s$ is defective, and its defectivity is $2$ in the first
case and $1$ in the second.}
\par
\bigskip
The proof of the Theorem uses the same kind of procedures as Proposition 2.3, but 
it also makes use of the $Horace$ $differential$ $method$ (see \ref
{AH}), and its proof is quite long since a lot of different cases
have to be considered. A complete proof can be found in \ref
{CGG3}.
\par
\bigskip
\par
\Coroll {2.6} {\it Let $a_1,a_2\geq 1$, and $V=V_{a_1,a_2}$ be the Segre-Veronese
$(a_1,a_2)$-embedding of $\PP 1\times \PP 1$. Then  $V$ is not Grassman defective, except when
$(a_1,a_2)=(2\alpha,1)$, and in this case $V$ is $(1,2\alpha )$-defective.}
\par
\medskip
\Proof\ By  Proposition 1.0, this statement is equivalent to
Theorem 2.5.\par \hfill  \qed
\par
\bigskip
\bigskip
\noindent {\bf 3. Other results.}
\par
\medskip
The correspondence between the dimension of $V^s_{\bf n,a}$ and the Hilbert function of a scheme 
$Z$ made of $s$ generic 2-fat points in $\PP {\bf n}$ 
(see Corollary 1.2) allows us to deduce results on $\dim V^s_{\bf n,a}$ from previous results on Segre Varieties
and from
properties of Hilbert funcions. Namely we have (notations as in \S 1):
\par
\bigskip
\noindent {\bf Proposition 3.1:} {\it Let $V_{\bf n,a}$ be a Segre-Veronese variety and $s\geq 2$ 
be such that $\dim V^s_{\bf n,a}= ns+s-1$ (the expected dimension). Then also
$V^s_{\bf n,b}$ has the (same) expected dimension for any ${\bf b}$ such that $b_i\geq a_i$ $\forall i=1,...,t$.}
\par
\medskip
$Proof:$  This is an immediate consequence of Corollary 1.2, since our hypothesis on $V^s_{\bf n,a}$ 
amounts to say that for a 
generic scheme $Z\subset \PP {n_1}\times...\times \PP {n_t}$ made of $s$ 2-fat points, 
$H(Z,{\bf a})= s(n+1)= lenght\ Z$
and this trivially implies that  $H(Z,{\bf b})= s(n+1)= lenght\ Z$ for all ${\bf b}$ such that $b_i\geq a_i$ $\forall
i=1,...,t$. 
\par
\medskip
From this follows:
\par
\medskip

\noindent {\bf Proposition 3.2:} {\it 
Let $1\leq n_1\leq n_2 \leq ... \leq n_t$, 
 and $V=V_{\bf n,a}$ be the Segre-Veronese
embedding of $\PP {n_1}\times ...\PP {n_t}$. If we are not in the case 
$t=2$, ${\bf a}=(1,1)$, then $\dim V^s = s(n_1+...+n_t+1)-1$ for all $s \leq n_1+1$.}
\par
\medskip
$Proof:$ 
If $t\geq 3$, from \ref {CGG2}, Proposition 2.3, we have that $H(Z,{\bf 1})=s(n_1+...+n_t+1)=lenght \ Z$, so we can
conclude by Proposition 3.1.
\par
When $t=2$, we know that for ${\bf a}=(1,1)$, $V^s$ is defective for all $2\leq s\leq n_1$ (e.g. see \ref {CGG2},
Proposition 2.3 again). We will be done if we show that (again by Proposition 3.1),   $V^{n_1+1}_{1,2}$ and
$V^{n_1+1}_{2,1}$ are not defective. 
\par
Without loss of generality, we can consider $Z$ such that $Supp\ Z = \{P_0,...,P_{n_1}\}\subset \PP {n_1}\times \PP {n_2}$, where $P_i$ is the coordinate point associated to the bihomogeneous ideal ${\goth p}_i = (x_0,...,\hat x_i,...,x_{n_1};y_0,...,\hat y_i,...,y_{n_2})$ in the ring $k[(x_0,...x_{n_1};y_0,...,y_{n_2}]$. Let $I={\goth p}_0^2\cap ...\cap{\goth p}_{n_1}^2$ be the ideal associated to $Z$; since $I$ is a monomial ideal (for the $P_i$'s are coordinate points), we only need to show that the monomials not in $I_{(1,2)}$ and those not in $I_{(2,1)}$ are $(n_1+1)(n_1+n_2+1)$ in number.
\par
The monomials not in $I_{(1,2)}$ are of two types:
\par
\noindent $x_iy_iy_j$, with $i=0,...,n_1$, $j=0,...,n_2$ ($(n_1+1)(n_2+1)$ of them);
\par
\noindent $x_iy_j^2$, with $i\neq j$, both in $\{0,...,n_1\}$ ($(n_1+1)n_1$ of them).
\par
For a total number of $(n_1+1)(n_1+n_2+1)$, as requested.
\par
The monomials not in $I_{(2,1)}$ are of two types:
\par
\noindent $x_ix_jy_k$, with $i\neq j$, both in $\{0,...,n_1\}$ and $k=i$ or $k=j$ ($(n_1+1)n_1$ of them);
\par
\noindent $x_i^2y_j$, 
with $i=0,...,n_1$, $j=0,...,n_2$ ($(n_1+1)(n_2+1)$ of them).
\par
For a total number of $(n_1+1)(n_1+n_2+1)$, as requested. \par \hfill  \qed
\par 
\bigskip
From this result we get immediately the following:
\par
\Coroll {3.3} {\it Let $r,k\geq 1$, $d\geq 2$, and $V$ be the d-ple
(Veronese) embedding of $\PP r$. Then $V$ is not (Grassman)
$(k,s-1)$-defective, for $s\leq min \{r+1,k+1\}$.}
\par
\medskip
\Proof\ By Proposition 1.0, this statement is equivalent to say that $\PP r\times \PP k$ in the Segre-Veronese embedding
of bidegree $(d,1)$ is not $(s-1)$-defective, hence the results follows from Proposition 3.2. \par \hfill  \qed
\par
\bigskip
\bigskip
Up to this point in this section we only proved results about Segre-Veronese varieties 
which are NOT defective. What follows is a list 
of examples of defective varieties; the way one can check the defectivity in all this examples is the same: we should
have $ V^s_{\bf n,a} = \PP N$, but instead it is easy to find a way to split ${\bf a} = {\bf b}+{\bf
c}=(b_1,...,b_t)+(c_1,...,c_t)$ in such a way that there is a form $f_1$ of multidegree $(b_1,...,b_t)$ and an $f_2$ of
multidegree $(c_1,...,c_t)$  passing through $s$ generic (simple) points, hence there is at least a form of degree ${\bf
a}$ (namely, $f_1f_2$) through $s$ generic 2-fat points which was not supposed to exist.
\par
\medskip
In the following list we always have $m\geq 1$, and we give values $s,{\bf n,a}$ for which $ V^s_{\bf n,a} $ is defective:
\par
\medskip
\noindent $\PP 1\times \PP m$, ${\bf a}=(2k,2)$, ${\bf b}={\bf c}=(k,1)$, $k\geq 1$, $s=\lceil{(2k+1)(m+1)\over
2}\rceil$;

\noindent $\PP 2\times \PP 2$, ${\bf a}=(2,2)$, ${\bf c}={\bf d}=(1,1)$, $s=8$;

\noindent $\PP 1\times \PP 1\times \PP m$, ${\bf a}=(1,1,2)$, ${\bf b}=(1,0,1)$, ${\bf c}=(0,1,1)$, $k\geq 1$, $s=2m+1$;

\noindent $\PP 1\times \PP m\times \PP m$, ${\bf a}=(2k,1,1)$, ${\bf b}=(k,1,0)$, ${\bf c}=(k,0,1)$, $k\geq 1$,
$s=km+k+m$;

\noindent $\PP 1\times \PP r\times \PP m$, ${\bf a}=(r+m,1,1)$, ${\bf b}=(r,1,0)$, ${\bf c}=(m,0,1)$, $s=rm+r+m$;

\noindent $\PP 1\times \PP 1\times \PP m$, ${\bf a}=(2,2,2)$, ${\bf b}={\bf c}=(1,1,1)$, $m\leq 3$, $s=4m+3$;

\noindent $\PP 2\times \PP m\times \PP m$, ${\bf a}=(2,1,1)$, ${\bf b}=(1,1,0)$, ${\bf c}=(1,0,1)$, $s=3m+2$;

\noindent $\PP 1\times \PP 1\times \PP 2\times \PP 5$, ${\bf a}=(2,1,1,1)$, ${\bf b}=(1,1,1,0)$, ${\bf
c}=(1,0,0,1)$, $s=11$;

\noindent $\PP 1\times \PP 1\times \PP 1\times \PP {2m-1}$, ${\bf a}=(m,1,1,1)$, ${\bf b}=(m-1,1,1,0)$, ${\bf
c}=(1,0,0,1)$, $m>1$, $s=4m-1$;

\noindent $\PP 1\times \PP 1\times \PP 1\times \PP {2m}$, ${\bf a}=(m,1,1,1)$, ${\bf b}=(m-1,1,1,0)$, ${\bf
c}=(1,0,0,1)$, $m\geq 4$, $s=4m-1$.

\medskip

Notice that the penultimate case for $m=1$ is defective too, but it is not of the same kind of all the others ($V^s$ has
dimension 13 and not 14); for this example see \ref {CGG 4}, Example 2.2.
\par
\bigskip
Of course from these examples we can derive examples of Grassmann defectivity, again by using Proposition 1.0; we will
just notice what we get from the last two cases:
\par
\medskip
\Coroll {3.4} {\it Let $V$ be the $(m,1,1)$-ple (Segre-Veronese) embedding of $\PP 1\times \PP 1\times \PP 1$. Then $V$
is (Grassman)
$(2m-1,4m-2)$-defective, for $m\geq 1$ and $(2m,4m-2)$-defective, for $m\geq 4$.}
\par
\medskip
The Corollary shows that the Segre Veronese varieties given by a $(2\alpha,1,1)$-embedding of $\PP 1\times \PP 1\times
\PP 1$ are both defective (see Theorem 2.5) and Grassmann defective.
\par
\bigskip

\noindent {\bf 4. Partially symmetric tensors.}
\par
\bigskip
\noindent Now we want to describe how to interpret our Segre-Veronese
embeddings from the point of view of tensors (e.g. see \ref {Ge} or \ref {IK} for the Veronese case and \ref {CGG2} 
for the Segre case). In
order to look at the Veronese variety (or $t$-uple embedding of
$\PP r$) $\nu_t (\PP n)\subset \PP N$, $N={t+r\choose r}$ one can
consider the Segre Variety $\nu (\PP r\times...\times \PP r)$, 
with $t$ factors, in $\PP M$, $M=(r+1)^t-1$, and then consider the
action of the symmetric group $S_t$ on $\PP M$ where, if the
variables in $\PP M$ are
$\{z_{(1,0,...,0)...,(1,0,...,0)},...,z_{(0,...,0,1),...(0,...,0,1)}\}$,
the action of an element $\sigma\in S_t$ is defined by $$ \sigma
(z_{(1,0,...,0)...,(1,0,...,0)},...,z_{(0,...,0,1),...(0,...,0,1)})=
(z_{\sigma ((1,0,...,0)...,(1,0,...,0))},...,z_{\sigma
((0,...,0,1),...(0,...,0,1))}). $$ The invariant subspace of $\PP
M$ with respect to this action is actually a linear space $\cong
\PP N$, and the linear equations which define it give the
required symmetries for the tensors parameterized by the points of
$\PP M$.  We can view $ \PP M$ as the parameterizing space of all
the $(r+1)^t$ tensors, and $\PP N$ inside it as the subspace of
symmetric ones: then the Segre Variety and the Veronese
parameterize the rank one (decomposable) tensors.
\par
Notice that the symmetric tensors of rank one correspond to forms
that can be written as powers of linear forms. Notice also that
when we say that a symmetric tensor has rank one, i.e that it is
decomposable, we mean that it is decomposable as an element of the
Tensor Algebra $V\otimes ... \otimes V=V^{\otimes t}$ (where $\PP
 r = \PPP (V)$), not of the symmetric algebra Sym$_t(V)$.
\par
Consider for example a rational normal curve $C_t\i \PP t$; we are
used to view its ideal as generated by the $2\times 2$ minors of a
$2\times t$ catalecticant matrix of indeterminates (or also by the
$2\times 2$ minors of a different catalecticant matrix, see e.g.
\ref {Pu}). From the point of view above we should look at the
ideal of the Segre embedding $V_t$: $(\PP 1)^t \rightarrow \PP
{2^t-1}$, which is generated by the $2\times 2$ minors of a
$2\times 2\times ...\times 2$ ($t$ times) tensor (e.g. see \ref
{Gr} and \ref {Ha}); the ideal of $C_t$ comes from the ideal of
$V_t$ modulo the symmetry relations (given by the action of the
symmetric group $S_t$ on $\PP {2^t-1}$) which define a linear
space $\PP t$ in $\PP {2^t-1}$.
\par


This can be thought as a more ``complete'' way to view those
ideals, with respect to the usual way (as given by minors of
catalecticant matrices) since the tensor represents ``more
faithfully'' their symmetries.
\par
Now consider e.g. the case $t=3$; we can think of ``stopping halfway''
between the Segre variety $V_3$ (parameterizing $2\times 2\times
2$ decomposable tensors in $\PP 7$) and the rational normal curve
$C_t$ (which parameterizes decomposable $2\times 2\times 2$
symmetric tensors) by considering the Segre-Veronese embedding
$V_{(2,1)}$ of $\PP 1\times \PP 1$ of degree $(2,1)$ into $\PP 5$.
\par
We can consider the action of the symmetric group $S_2$ on $\PP 7$
which symmetrizes its variables $x_{ijk}$, $i,j,k\in {0,1}$ only
with respect to $i$ and $j$. The invariant space for this action
is a linear space $\PP 5\i \PP 7$, and it cuts $V_3$ exactly in
$V_{(2,1)}$. Hence the variety $V_{(2,1)}$ parameterizes $2\times
2\times 2$ ``partially symmetric tensors'', i.e. tensors whose
entries are symmetric only with respect to the first two indeces.
\par
\medskip
In general, consider $(\PP r)^t=\PP r\times ...\times \PP r$, $t$
times, and its Segre-Veronese embedding $V_{{\bf r},{\bf a}}$, ${\bf r}=(r,...,r)$ and 
${\bf a}=(a_1,...,a_t)$, into the space $\PP N$. Let $a
= a_1+...+a_t$, and consider the Segre embedding of $(\PP r)^a$
into $\PP M$, where $M=(r+1)^a-1$. We can view $\PP N$ inside $\PP
M$ as the space of tensors which are invariant with respect to the
actions of $S_{a_1}$,...,$S_{a_t}$ on the variables of relative
indeces. So those are ``partially symmetric" tensors (for $t=1$ we
get symmetric tensors and $V_{{\bf r},{\bf a}}$ is the Veronese
variety, while for $a_1=...=a_t=1$ they are generic tensors and
$V_{{\bf r},{\bf a}}$ is the Segre embedding of $(\PP r)^t$).
\par
So the Segre-Veronese variety $V_{{\bf r},{\bf a}}$, will
parameterize the partially symmetric tensors (with respect to the
actions of $S_{a_1}$,...,$S_{a_s}$) in $\PP M$ which are
decomposable. Since those are the tensors of tensor rank 1 (e.g.
see \ref {CGG2}), the secant varieties of  $V_{{\bf r},{\bf a}}$
give the stratification by tensor rank of those partially
symmetric tensors.
\par
\medskip
\bigskip
\bigskip

\centerline {{\bf REFERENCES}}

\medskip\noindent[{\bf AH} ]: J. Alexander, A. Hirschowitz. {\it An asymptotic vanishing theorem for generic unions of multiple points.} Inv. Math. {\bf 140} (2000), 303-325.
\par
\medskip
\noindent [{\bf BCS}]: P. B\"{u}rgisser, M. Clausen, M.A.
Shokrollahi, {\it Algebraic Complexity Theory}, Vol. 315, Grund.
der Math. Wiss., Springer, 1997
\par
\medskip
\noindent [{\bf BuM}] B\u arc\u anescu, \c S.; Manolache, {\it N. Betti numbers of Segre-Veronese singularities.}  Rev. Roumaine Math. Pures Appl.  {\bf 26}  (1981), 549--565. 
\par
\medskip

\noindent [{\bf CEG}]: M.V.Catalisano, P.Ellia, A.Gimigliano.
 {\it Fat points on rational normal curves.}  J. of Algebra,  {\bf 216},
(1999), 600-619.
\par
\medskip

\noindent [{\bf CGG1}]: M.V.Catalisano, A.V.Geramita,
A.Gimigliano. {\it On the Secant Varieties to the Tangential
Varieties of a Veronesean.} Proc. A.M.S. {\bf 130} (2001).
975-985.
\par
\medskip
\noindent [{\bf CGG2}]: M.V.Catalisano, A.V.Geramita,
A.Gimigliano. {\it Rank of Tensors, Secant Varieties of Segre
Varieties and Fat Points.}  Linear Alg. Appl. {\bf 355}, (2002), 261-285.

\noindent [{\bf CGG3}]: M.V.Catalisano, A.V.Geramita,
A.Gimigliano. {\it Higher Secant varieties of Segre embeddings of $\PP 1\times\PP 1\times \PP 1$.} Preprint (2003). 

\noindent [{\bf CGG4}]: M.V.Catalisano, A.V.Geramita,
A.Gimigliano. {\it Higher Secant varieties of Segre varieties of $\PP 1\times ...\times \PP 1$.} Preprint (2003). 

\medskip\noindent {\bf COCOA}: A. Capani, G. Niesi, L. Robbiano,
{\it CoCoA, a system for doing Computations in Commutative
Algebra} (Available via anonymous ftp from: cocoa.dima.unige.it).

\medskip
\noindent [{\bf Cha}]: K.Chandler. {\it  A brief proof of a
maximal rank Theorem for generic double points in projective
space} Trans. Amer. Math. Soc. {\bf 353} (2000), 1907-1920

\noindent [{\bf ChCi}]: L.Chiantini, C.Ciliberto. {\it Weakly
defective varieties.} Trans. Am. Math. Soc. {\bf 354} (2001).
151-178

\medskip
\noindent [{\bf ChCo}]: L.Chiantini, M.Coppens. {\it Grassmannians
for secant varieties.} Forum Math. {\bf 13} (2001) 615-628.
\par
\medskip

\noindent [{\bf DF}]: C. Dionisi, C.Fontanari. {\it Grassmann
defectivity \`a la Terracini.} Preprint (AG-0112149).
\par
\medskip

\noindent [{\bf E}]: R.Ehrenborg. {\it On Apolarity and Generic
Canonical Forms.} J. of Algebra {\bf 213} (1999), 167-194.

\medskip\noindent [{\bf ER}]: R.Ehrenborg, G.-C. Rota. {\it Apolarity and
canonical forms for homogeneous polynomials.} European J. of
Combinatorics {\bf 14} (1993), 157-181.

\medskip\noindent [{\bf GHKM}]: D.Geiger, D.Hackerman, H.King, C.Meek.
{\it Stratified Exponential Families: Graphical Models and Model
Selection.} Annals of Statistics, {\bf 29} (2001), 505-527.

\medskip\noindent [{\bf Ge}]: A.V.Geramita.
{\it Inverse Systems of Fat Points}, Queen's Papers in Pure and
Applied Math.  {\bf 102},  {\it The Curves Seminar at Queens',
vol. X} (1998).

\medskip\noindent [{\bf Ge2}]: A.V.Geramita. {\it Catalecticant Varieties.},
Lecture Notes in Pure and Applied Math. Dekker {\bf 206}.

\medskip\noindent [{\bf Hr}]: J.Harris. {\it Algebraic Geometry, a
First Course.} Springer-Verlag, New York (1993).

\medskip\noindent [{\bf IK}]: A.Iarrobino, V.Kanev. {\it Power Sums,
Gorenstein algebras, and determinantal loci.} Lecture Notes in
Math. {\bf 1721}, Springer, Berlin, (1999).

\medskip\noindent [{\bf Id}]: M.Id\`a. {\it Private communication.}

\medskip\noindent [{\bf La}]:A. Laface, {\it On linear systems of curves on rational scrolls}, Geom.
Dedicata {\bf 90} (2002), 127-144.

\medskip\noindent [{\bf Pa}]: F.Palatini. {\it Sulle variet\`a algebriche
per le quali sono di dimensione minore} {\it dell' ordinario,
senza riempire lo spazio ambiente, una o alcuna delle variet\`a}
{\it formate da spazi seganti.} Atti Accad. Torino Cl. Scienze
Mat. Fis. Nat. {\bf 44} (1909), 362-375.

\medskip
\noindent [{\bf Te}]: A.Terracini. {\it Sulle} $V_k$ {\it per cui
la variet\`a degli} $S_h$ $(h+1)${\it -seganti ha dimensione
minore dell'ordinario.} Rend. Circ. Mat. Palermo {\bf 31} (1911),
392-396. par
\medskip
\noindent [{\bf Te2}]: A.Terracini. {\it Sulla rappresentazione
delle coppie di forme ternarie mediante somme di potenze di forme
lineari.} Ann. Mat. Pur ed appl. {\bf XXIV}, {\bf III} (1915),
91-100.

\medskip\noindent [{\bf Z}]: F.L.Zak. {\it Tangents and Secants of Algebraic
Varieties.} Translations of Math. Monographs, vol. {\bf 127} AMS.
Providence (1993).

\bigskip

{\it M.V.Catalisano, Dip. Matematica, Univ. di Genova, Italy.}

{\it e-mail: catalisa@dima.unige.it}

\medskip
{\it A.V.Geramita, Dept. Math. and Stats. Queens' Univ. Kingston,
Canada} {\it and Dip. di Matematica, Univ. di Genova. Italy.}

{\it e-mail: geramita@dima.unige.it ; tony@mast.queensu.ca}

\medskip
{\it A.Gimigliano, Dip. di Matematica and C.I.R.A.M., Univ. di
Bologna, Italy.}

{\it e-mail: gimiglia@dm.unibo.it}

\end